\documentclass[11pt]{article}
\usepackage{graphicx, amsmath, amsfonts, amssymb, amsthm}
\newtheorem{theorem}{Theorem}[section]
\newtheorem{corollary}[theorem]{Corollary}
\newtheorem{lemma}[theorem]{Lemma}

\theoremstyle{definition}
\newtheorem{definition}[theorem]{Definition}

\DeclareMathOperator{\id}{id}
\DeclareMathOperator{\Graph}{Graph}

\begin{document}

%%%%%%%%%%%%%%%%%%%%%%%%%%%%%%%%%%%%%%%%%%%%%%%%%%%
%%%%%                                         %%%%%
%%%%%     On the first integral conjecture    %%%%%
%%%%%            of René Thom                 %%%%%
%%%%%                                         %%%%%
%%%%%  J. Cresson, A. Daniilidis, M Shiota    %%%%%
%%%%%                                         %%%%%
%%%%%            October 8, 2007              %%%%%
%%%%%          (version preprint)             %%%%%
%%%%%%%%%%%%%%%%%%%%%%%%%%%%%%%%%%%%%%%%%%%%%%%%%%%

%\usepackage[notref,notcite]{showkeys}

\begin{center}
{\LARGE On the First Integral Conjecture of Ren\'{e} Thom}

\vspace{1cm}

{\large \textsc{Jacky\ CRESSON, Aris DANIILIDIS \& Masahiro SHIOTA}}
\end{center}

\bigskip

\noindent\textbf{Abstract.} More than half a century ago R.~Thom
asserted in an unpublished manuscript that, generically, vector
fields on compact connected smooth manifolds without boundary can
admit only trivial continuous first integrals. Though somehow
unprecise for what concerns the interpretation of the word
\textquotedblleft generically\textquotedblright, this statement is
ostensibly true and is nowadays commonly accepted. On the other
hand, the (few) known formal proofs of Thom's conjecture are all
relying to the classical Sard theorem and are thus requiring the
technical assumption that first integrals should be of class~$C^{k}$
with $k\geq d,$ where $d$ is the dimension of the manifold. In this
work, using a recent nonsmooth extension of Sard theorem we
establish the validity of Thom's conjecture for locally Lipschitz
first integrals, interpreting genericity in the $C^{1}$ sense.

\bigskip

\noindent\textbf{Key words. }Structural stability, first integral, o-minimal
structure, Sard theorem.

\bigskip

\noindent\textbf{AMS Subject Classification.} \ \textit{Primary} 37C20 ;
\textit{Secondary} 34D30, 14P10.

\section{Introduction}

The purpose of this paper is to discuss the following conjecture
attributed to Ren\'{e} Thom (see~\cite{Thom}, \cite{Francoise1995},
\cite{Peixoto1967}, for example) which is part of the
\textit{folklore} in the dynamical systems community:

\bigskip

\noindent\textbf{Thom conjecture}: \textit{For $1\leq r\leq\infty$,
$C^{r} $-generically vector fields on }$d$\textit{-dimensional
compact, smooth, connected manifolds without boundary do not admit
nontrivial continuous first integrals.}

\bigskip

Ren\'{e} Thom \cite{Thom} proposes a scheme for a formal proof
relying on the assumption that a $C^{r}$~closing lemma ($r\geq1$) is
true \cite{Francoise1995}. The $C^{1}$-closing lemma (case~$r=1$)
has indeed been proved by Pugh (\cite{Pugh1967a}, \cite{Pugh1967b}).
Nevertheless, very little is known for a $C^r$-closing lemma with
$r\geq 2$ (see also \cite{Gutierrez1987}, \cite{Pugh1975}) so that
Thom's strategy should be revised.

\medskip

In \cite{Peixoto1967}, Peixoto proves Thom conjecture for $r=1$ assuming that
first integrals are of class $C^{k}$ with $k\geq d$, \textit{i.e.} a relation
between the regularity class of the first integrals to be considered and the
dimension of the underlying compact manifold. As pointed out by Peixoto, this
condition is only of technical nature and relates to the use of the classical
Sard theorem in a crucial part of the proof. More precisely, denoting by
$\mathcal{X}_{M}$ the set of $C^{1}$-vector fields on the $d$-dimensional
compact manifold~$M$, Peixoto's proof is divided in three steps:

\bigskip

Let $X\in\mathcal{X}_{M}$ and $f\colon M\rightarrow{\mathbb{R}}$ be a first integral
of class $C^{d}$.

\medskip

- By Sard's lemma, there exists in $f(M)$ an interval $]a,b[$ made
up of regular values. For any~$y\in]a,b[$, $f^{-1}(y)$ is an
$(d-1)$-dimensional, compact, differentiable manifold, invariant
under~$X$.

\medskip

- By Pugh's general density theorem for each $y\in]a,b[$, $f^{-1}(y)$ does not
contain singularities or closed orbits of $X$ since they are generic. As a
consequence, singularities and closed orbits are all located at the critical
levels of $f$.

\medskip

- Any trajectory $\gamma$ in $f^{-1}(y)$ is such that $\omega(\gamma)\subset
f^{-1}(y)$ and cannot be contained in the closure of the set of singularities
or closed orbits, in contradiction to Pugh's general density theorem.

\bigskip

T.~Bewley \cite{Bewley1971} extends Peixoto's theorem for $1\leq r\leq\infty$.
The proof has then been simplified by R.~Ma\~{n}e \cite{Mane1973}. Ma\~{n}e's
proof seems to have been rediscovered by M.~Hurley \cite{Hurley1986}. However,
the technical assumption of Peixoto ($C^{d}$-regularity of the first
integrals) stays behind all these works, because of the use of Sard's theorem.
Nevertheless, Thom conjecture seems to be true in general.

\medskip

In this paper, we cancel the aforementioned regularity condition by
interpreting the word ``nontrivial'' as ``being essentially
definable with respect to an o-minimal approximation'' (see
Definition~\ref{Definition_essential} below). In this context we
prove the validity of Thom's conjecture for $r=1$ and for Lipschitz
continuous first integrals.

\smallskip

The technic of the proof follows the strategy used by M.~Artin and
B.~Mazur (\cite{Artin-Mazur1965}) to prove that generically the
number of isolated periodic points of a diffeomorphism grows at most
exponentially. Indeed, using a recently established nonsmooth
version of Sard theorem (see \cite[Theorem~7]{BDLS2005} or
\cite[Corollary~9]{BDLS2007}) and Peixoto's scheme of proof
(\cite{Peixoto1967}) we first show that in an o-minimal manifold
(that is, a manifold that is an o-minimal set), generically,
o-minimal first integrals are constant. Then by approximating every
compact differentiable manifold by a Nash manifold we derive a
general statement.

\bigskip

\begin{itemize}
\item \textbf{Preliminaries in dynamical systems.}
\end{itemize}

\smallskip

Given a $C^{1}$ manifold $M$ we denote by $\mathcal{X}_{M}$ the
space of all $C^{1}$-vector fields on $M$ equipped with the $C^{1}$
topology. Let $\phi _{t}\colon M\rightarrow M$ be the one-parameter group
of diffeomorphisms generated by a vector field $X$ on $M$. A point
$p\in M$ is called \textit{nonwandering}, if given any neighborhood
$U$ of~$p$, there are arbitrarily large values of $t$ for which
$U\cap\phi_{t}(U)\not =\emptyset$. Denoting by $\Omega$ the set of
all nonwandering points we have the following genericity result due
to C.~Pugh (\cite{Pugh1967a}, \cite{Pugh1967b}).

\medskip

\noindent
\textbf{General density theorem (GDT):} The set $\mathcal{G}_{M}$ of
vector fields $X\in\mathcal{X}_{M}$ such that properties
(G$_{1}$)--(G$_{4}$) below hold is residual in $\mathcal{X}_{M}$.

\smallskip

(G$_{1}$) $X$ has only a finite number of singularities, all generic
;\smallskip

(G$_{2}$) Closed orbits of $X$ are generic ;\smallskip

(G$_{3}$) The stable and unstable manifolds associated to the singularities
and the closed orbits of $X$ are transversal ;\smallskip

(G$_{4}$) $\Omega=\bar{\Gamma}$, where $\Gamma$ stands for the union of all
singular points and closed orbits of $X$.

\bigskip

We use the following definition of a \textit{first integral}:

\begin{definition} [First integral] A first integral of a vector field $X$ on a
compact connected manifold $M$ of dimension $d$ is a continuous
function $f\colon M\rightarrow\mathbb{R}$ which is constant on the orbits
of the flow generated by $X$ but it is not constant on any nonempty
open set of $M$.
\end{definition}

As mentioned in the introduction, Peixoto \cite{Peixoto1967} only considers
$C^{k}$-first integrals with $k\geq d$.

\medskip
%\newpage
\begin{itemize}
\item \textbf{Preliminaries in o-minimal geometry.}
\end{itemize}

Let us recall the definition of an o-minimal structure (see
\cite{Dries-Miller96} for example).

\begin{definition}
[o-minimal structure]\label{Definition_o-minimal} An o-minimal
structure on the ordered field $\mathbb{R}$ is a sequence of Boolean
algebras $\mathcal{O}=\{\mathcal{O}_{n}\}_{n\geq1}$ such that for
each $n\in \mathbb{N\smallskip}$

\begin{enumerate}
\item[(i)] $A\in\mathcal{O}_{n}$ $\Longrightarrow$ $A\times\mathbb{R}
\in\mathcal{O}_{n+1}$ and $\mathbb{R}\times A\in\mathcal{O}_{n+1}$
$;$

\item[(ii)] $A\in\mathcal{O}_{n+1}\Longrightarrow$ $\Pi(A)\in
\mathcal{O}_{n}$
\newline$(\Pi\colon \mathbb{R}^{n+1}\rightarrow\mathbb{R}^{n}$ denotes the
canonical projection onto $\mathbb{R}^{n})$ ;

\item[(iii)] $\mathcal{O}_{n}$ contains the family of algebraic
subsets of $\mathbb{R}^{n}$, that is, the sets of the form
\[
\{x\in\mathbb{R}^{n}:p(x)=0\},
\]
where $p\colon \mathbb{R}^{n}\rightarrow\mathbb{R}$ is a polynomial
function $;$

\item[(iv)] $\mathcal{O}_{1}$ consists exactly of the finite unions of
intervals and points.
\end{enumerate}
\end{definition}

An important example of o-minimal structure is the collection of
\emph{semialgebraic sets} (see \cite{Coste99} for example), that is,
sets that can be obtained by Boolean combinations of sets of the
form
\[
\left\{  x\in\mathbb{R}^{n}:\,p(x)=0,\,q_{1}(x)<0,\ldots,\,q_{m}(x)<0\right\}
,
\]
where $p,q_{1},\ldots,q_{m}$ are polynomial functions in $\mathbb{R}^{n}$.
Indeed, properties (i), (iii) and (iv) of Definition~\ref{Definition_o-minimal}
are straightforward, while (ii) is a consequence of the Tarski-Seidenberg principle.

\smallskip

A subset $A$ of $\mathbb{R}^{n}$ is called \emph{definable} (in the o-minimal
structure $\mathcal{O}$) if it belongs to $\mathcal{O}_{n}.$ Given any
$S\subset\mathbb{R}^{n}$ a mapping $F\colon S\rightarrow\mathbb{R}$ is called
\emph{definable} in $\mathcal{O}$ (respectively, semialgebraic) if its graph
is a definable (respectively, semialgebraic) subset of~$\mathbb{R}^{n}
\times\mathbb{R}$.

\medskip

\begin{itemize}
\item \textbf{Preliminaries in variational analysis.}
\end{itemize}

Let $g\colon U\rightarrow\mathbb{R}$ be a Lipschitz continuous function
where $U$ is a nonempty open subset of~$\mathbb{R}^{d}$. The
\emph{generalized derivative} of $g$ at $x_{0}$ in the direction
$v\in\mathbb{R}^{n}$ is defined as follows (see
\cite[Section~2]{Clarke1983} for example):
\begin{equation}
g^{o}(x_{0},e)=\,\limsup_{\begin{subarray}{c}
x\rightarrow x_{0}\\
t\searrow0^{+}\end{subarray}}
\frac{g(x+te)-g(x)}{t} \label{Clarke derivative}
\end{equation}
where $t\searrow0^{+}$ indicates the fact that $t>0$ and
$t\rightarrow0$. It turns out that the function $v\mapsto
g^{o}(x_{0},v)$ is positively homogeneous and convex, giving rise to
the \emph{Clarke subdifferential} of $g$ at $x_{0}$ defined as
follows:
\begin{equation}
\partial g(x_{0})\,=\,\{\,x^{\ast}\in\mathbb{R}^{d}:\,g^{o}(x_{0}
,v)\geq\langle x^{\ast},v\rangle,\,\forall
v\in\mathbb{R}^{d}\}\text{.} \label{subdiff}
\end{equation}
In case that $g$ is of class $C^{1}$ (or more generally, strictly
differentiable at $x_{0}$) it follows that $$\partial
g(x_{0})=\{\nabla g(x_{0})\}.$$

\smallskip

A point $x_{0}\in U$ is called \emph{Clarke critical}, if
$0\in\partial g(x_{0}).$ We say that $y_{0}\in g(U)$ is a
\emph{Clarke critical value} if the level set $g^{-1}(y_{0})$
contains at least one Clarke critical point. Given a Lipschitz
continuous function $f\colon M\rightarrow\mathbb{R}$ defined on a $C^{1}$
manifold $M$ we give the following definition of (nonsmooth)
critical value.

\begin{definition}
[Clarke critical value]\label{Definition_critical} We say that
$y_{0}\in f(M)$ is a Clarke critical value of the function
$f\colon M\rightarrow\mathbb{R}$, if there exists $p\in f^{-1}(y_{0})$ and
a local chart $(\varphi,U)$ around $p$ such that
$0\in\partial(f\circ\varphi^{-1})(\varphi(p))$. In this case, $p\in
M$ is a Clarke critical point for $f$.
\end{definition}

\noindent It can be easily shown (see \cite[Exercise 10.7]{RW1998},
for example) that the above definition does not depend on the choice
of the chart.

\section{Main results}

Throughout this section $M$ will be a $C^{1}$ compact connected
submanifold of $\mathbb{R}^{n}$ (without boundary), $TM$ its
corresponding tangent bundle and $\mathcal{X}_{M}$ the space of
$C^{1}$ vector fields on $M$ equipped with the $C^{1}$ topology. Let
us recall that submanifolds of $\mathbb{R}^{n}$ admit $\varepsilon
$-tubular neighborhoods for all $\varepsilon>0$ sufficiently small.
We further consider a $C^{1}$-submanifold $N$ of $\mathbb{R}^{n}$, a
$C^{1} $-diffeomorphism $F\colon M\rightarrow N$ and $\varepsilon>0.$

\begin{definition}
[Approximation of a manifold] $(i)$ We say that $(N,F)$ is a $C^{1}
$-appro\-xi\-ma\-tion of $M$ (of precision $\epsilon$), if $N$
belongs to an $\epsilon$-tubular neighborhood $U_{\epsilon}$ of $M$
and $F$ can be extended to a $C^{1}$-diffeomorphism
$\tilde{F}\colon \mathbb{R}^{n}\rightarrow\mathbb{R}^{n}$, which is
isotopic to the identity $\id$, satisfies
$\tilde{F}|_{\mathbb{R}^{n}\diagdown U_{r}}\equiv \id$ and
\[
\max_{x\in\mathbb{R}^{n}}
\,\left\{
\,||\tilde{F}(x)-x||\,+\,||\,d\tilde{F}(x)-\id||\,\right\}
\,<\,\epsilon.
\]
(We shall use the notation $\tilde{F}\sim_{\epsilon}\id$ to
indicate that $\tilde{F}$ is $\epsilon$-$C^{1}$-closed to the
identity mapping.)\smallskip

$(ii)$ A $C^{1}$-approximation $(N,F)$ of $M$ is called
semialgebraic (respectively, definable) if the manifold $N$ is a
semialgebraic subset of $\mathbb{R}^{n}$ (respectively, a definable
set in an o-minimal structure).
\end{definition}

In the sequel, we shall need the following approximation result.

\begin{lemma}
[Semialgebraic approximation]\label{Lemma_shiota}Let $M$ be a $C^{1}$ compact
submanifold of $\mathbb{R}^{n}$. Then for every $\varepsilon>0$, there exists
a semialgebraic $\varepsilon$-approximation of $M$.
\end{lemma}

\noindent\textbf{Proof.} Fix $\varepsilon>0$ and let $U$ be an open
$\epsilon$-tubular neighborhood of $M$ in $\mathbb{R}^{n}$ for some
$\epsilon\in(0,\varepsilon)$. Applying
\cite[Theorem~I.3.6]{Shiota1987} (for $A=\mathbb{R}^{n}$ and
$B=C^{1}$), we deduce the existence of a $C^{1} $-embedding $F$ of
$M$ into $U$ which is $\epsilon$-close to the identity map $\id$ in
the $C^{1}$ topology such that $F(M)=N$ is a Nash manifold (that is,
$N$ is a $C^{\infty}$-manifold and a semialgebraic set). Then $F$~can be extended to a $C^{1}$ diffeomorphism $\tilde{F}$ of
$\mathbb{R}^{n}$ by a partition of unity of class $C^{1}$ such that
$\tilde{F}=\id$ on $\mathbb{R} ^{n}\setminus U$. Moreover there
exists a $C^{1}$ isotopy $\{F_{t} \}_{t\in\lbrack0,1]}$, such that
$F_{t}=\id$ on $\mathbb{R}^{n}\setminus U$, $F_{0}\equiv \id$ and
$F_{1}=\tilde{F}$ and the map $F_{t}\colon \mathbb{R}^{n}
\times\lbrack0,1]\rightarrow\mathbb{R}^{n}$ is $\epsilon$-close to
the projection to $\mathbb{R}^{n}$ in the $C^{1}$
topology.\hfill$\square$

\bigskip

Given a $C^{1}$-manifold $M$ and a $C^{1}$-vector field
$X\in\mathcal{X} _{M},$ the following result relates generic
singularities of hyperbolic type of $X$ with Clarke critical values
of Lipschitz continuous first integrals of $X$.

\begin{lemma}
[Location of singularities]\label{Lemma_critical} Assume
$f\colon M\rightarrow \mathbb{R}$ is a Lipschitz continuous first integral
for the vector field $X\in$ $\mathcal{X}_{M}$. Then all generic
singularities and all closed orbits of hyperbolic type are located
at the Clarke-critical level sets.
\end{lemma}

\noindent\textbf{Proof.} Let $p_{0}$ be either a singular point of
hyperbolic type or any point of a closed orbit of hyperbolic type
and let $\pi\colon U\rightarrow M$ be the exponential mapping around
$p_{0}=\pi(0),$ where $U$ is an open neighborhood of $0\in
T_{p_{0}}M\cong\mathbb{R}^{d}$ ($d$ denoting the dimension of $M$).
It follows that the function $g=f\circ\pi$ is Lipschitz continuous.
Since the stable and unstable manifolds of the flaw of the field $X$
at $p_{0}$ are transversal and since $f$ is a first integral, it
follows that for some basis $\{e_{i}\}_{i\in\{1,\ldots,d\}}$ of
$T_{p_{0}}M$ one has
\[
g^{o}(0,\pm e_{i})\geq g'(0,\pm e_{i}):=\lim_{t\searrow0+}
\frac{g(\pm te_{i})-g(0)}{t}=0,
\]
where $g^{0}(0,\cdot)$ is given by (\ref{Clarke derivative}). In
view of (\ref{subdiff}) we deduce that $0\in\partial g(0),$ hence
$f(p_{0})$ is a critical value of $f$.\hfill$\square$

\bigskip

\begin{lemma}
[Density of critical values for GDT fields]\label{Lemma_dense}In the
situation of Lemma~\ref{Lemma_critical}, let us further assume that
$X\in \mathcal{G}_{M}$. Then the Clarke critical values of $f$ are
dense in $f(M).$
\end{lemma}

\noindent\textbf{Proof.} Let $\Gamma$ denote the union of all
singular points and closed orbits of the field $X.$ Since
$X\in\mathcal{G}_{M}$, it follows from Lemma~\ref{Lemma_critical}
that the set $f(\Gamma)$ is included to the Clarke critical values.
Continuity of $f$ and compactness of $M$ yield
$\overline{f(\Gamma)}=f(\overline{\Gamma})=f(\Omega).$ Since
$\Omega$ contains all $\omega$-limits of orbits of $X,$ taking any
$y\in f(M)$ and any $x\in f^{-1}(y)\diagdown\Gamma$ we denote by
$\gamma$ the orbit passing through $x$ and by $\gamma_{\infty}$ the
set of $\omega$-limits of $\gamma$. Then by continuity
$y=f(\gamma)=f(\gamma_{\infty})\subset f(\Omega)$. This proves the
assertion.\hfill$\square$

\bigskip

\begin{corollary}
[Thom conjecture -- definable version]\label{Corollary_thom}Let
$X\in \mathcal{G}_{M}.$ Then $X$ does not admit any Lipschitz
continuous definable first integral.
\end{corollary}

\noindent\textbf{Proof.} Assume $f$ is a Lipschitz continuous first
integral of $X$ on $M$ and denote by $S$ the set of its Clarke
critical points. If $f$ is o-minimal, then so is $M$ (cf.~property
(ii) of Definition~\ref{Definition_o-minimal}), the tangential
mappings $\pi\colon U\subset T_{p_{0}}M\rightarrow M$ (around any point
$p_{0}\in M$) and the composite functions of the form $g=f\circ\pi$
(notation according to the proof of Lemma \ref{Lemma_critical}).
Note that $p$ is a critical point of $f$ if and only if
$\pi^{-1}(p)$ is a critical point of $g=f\circ\pi$ where $\pi$ is
any tangential mapping with $p\in\pi(U)$. Applying
\cite[Corollary~8]{BDLS2007} we deduce that the set of Clarke
critical values of each function $g$ is of measure zero. Using a
standard compactness argument we deduce that $f(S)$ is of measure
zero, thus in particular $f(M)\setminus f(S)$ contains an interval
$(y_{1},y_{2}).$ But this contradicts the density result of
Lemma~\ref{Lemma_dense}.\hfill$\square$

\bigskip

If the manifold $M$ is not a definable subset of $\mathbb{R}^{n}$ the above
result holds vacuously and gives no information. To deal with this case, the
forthcoming notion of essential o-minimality with respect to a given
o-minimal approximation turns out to be a useful substitute for our purposes.
Let us fix $\epsilon>0$ and a definable $\epsilon$-approximation $(N,F)$ of
$M$.

\begin{definition}
[Essential o-minimality with respect to a definable
approximation]\label{Definition_essential} A mapping
$f\colon M\rightarrow\mathbb{R}$ is called essentially o-minimal with
respect to a definable approximation~$(N,F)$ of~$M$ if the mapping
$f\circ F^{-1}\colon N\rightarrow\mathbb{R}$ is o-minimal.
\end{definition}

\noindent Note that every o-minimal function on $M$ is essentially
o-minimal with respect to any approximation $(N,F)$ of $M$ for which
the diffeomorphism $F$ is o-minimal. Setting
$$M=\{p\in\mathbb{R}^{2}:p\in\Graph(h)\}$$ where
$h(t)=t^{3}\sin(x^{-1}),$ if $t\neq0$ and $0$ if $t=0$, we obtain a
nondefinable $C^{1}$-submanifold of~$\mathbb{R}^{2}$. Thus the
(projection) function $f\colon M\rightarrow\mathbb{R}$ with $f((t,h(t))=t$
is not o-minimal. It can be easily seen that for every $\epsilon>0$
there exists an $\epsilon $-approximation $(N,F)$ with respect to
which $f$~is essentially o-minimal. On the other hand, if
$\chi_{K}$ is the characteristic function of the Cantor set $K$ of
$(0,1)$, the function $g\colon \mathbb{R}\rightarrow\mathbb{R}$ defined by
$g(x)=\int_{0}^{x}\chi_{K}(t)\,dt$ for all $x\in\mathbb{R}$ is not
essentially o-minimal with respect to any approximation. Roughly
speaking, a function that is not essentially o-minimal contains
intrinsic irreparable oscillations.

\bigskip

In view of Lemma~\ref{Lemma_shiota}, for every $\epsilon>0$ there
exists a $C^{1}$ definable manifold $N$ and a diffeomorphism
$F\colon M\rightarrow N$ such that $F\sim_{\epsilon}\id.$ Fixing the
approximation, we associate to every vector field $X\colon M\rightarrow
TM$ on $M$ the conjugate $C^{1}$-vector field
$\tilde{X}\colon N\rightarrow TN$ on $N$ defined as follows:
\[
\tilde{X}(q)=dF(F^{-1}(q),X(F^{-1}(q)).
\]
Note that $\tilde{X}$ is uniquely determined by $X$. Let us further denote by
$\mathcal{G}_{N}$ the vector fields of $N$ that satisfy the generic
GDT\ assumptions. We are ready to state our main result.

\begin{theorem}
[Genericity of non-existence of first integrals]Let $M$ be a $C^{1}$ compact
submanifold of $\mathbb{R}^{n}$ and $\epsilon>0.$ For the $C^{1}$ topology,
the set of vector fields in $M$ that do not admit Lipschitz continuous first
integrals which are essentially o-minimal with respect to a given definable
$\epsilon$-approximation of $M$ is generic.
\end{theorem}

\noindent\textbf{Proof.} Let us fix any definable
$\epsilon$-approximation of $M$ and let us denote by
$\mathcal{G}_{N}$ the vector fields of $N$ that satisfy the generic
GDT\ assumptions. By Pugh's density theorem $\mathcal{G} _{N}$ is a
$C^{1}$-residual subset of $\mathcal{X}_{N}$ and by
Corollary~\ref{Corollary_thom}, if $Y\in\mathcal{G}_{N}$ then $Y$
does not possess any o-minimal Lipschitz continuous first integral.
Let $\mathcal{G}$ denote the set of vector fields of $X$ that
conjugate inside $\mathcal{G}_{N} $, that is,
\[
\mathcal{G}=\{X\in\mathcal{X}_{M}:\;\tilde{X}\in\mathcal{G}_{N}\}.
\]
Then $\mathcal{G}$ is residual in $\mathcal{X}_{M}.$ Pick any $X\in\mathcal{G}$ and assume that
$f\colon M\rightarrow \mathbb{R}$ is a Lipschitz continuous first integral
of $X.$ Since the trajectories of $X$ are transported to the
trajectories of $\tilde{X} \in\mathcal{G}_{N}$ through the mapping
$F,$ it follows that $\tilde{f}=f\circ F$ is a first integral of
$\tilde{X}$ in $N.$ This shows that $f$ cannot be essentially
o-minimal with respect to $(N,F)$.\hfill$\square$

\bigskip
\noindent--------------------------------------
\medskip

\noindent\textbf{Acknowledgment.} Part of this work has been made
during a research visit of the second author to the University of
Pau (June 2007). The second author wishes to thank his hosts for
hospitality and Olivier Ruatta (University of Limoges) for useful
discussions. The first author thanks Jean-Pierre Bourguignon for his
help with R.~Thom's archive and J-P.~Fran\c{c}oise for useful
references. J.C. has been supported by the French projet de l'Agence Nationale de la Recherche "Int\'egrabilit\'e r\'eelle et complexe en M\'ecanique Hamiltonienne", N.JC05-41465

%\bigskip
\newpage

\noindent--------------------------------------

\noindent Jacky CRESSON

\smallskip

\noindent Laboratoire de Math\'{e}matiques Appliqu\'{e}es, UMR 5142 du CNRS\newline
Universit\'{e} de Pau et des Pays de l'Adour\newline avenue de
l'universit\'{e}, F-64000 Pau, France.\newline
and \newline
Institut de M\'ecanique C\'eleste et de Calcul des \'Eph\'em\'erides (IMCCE)\newline
UMR 8028 du CNRS - Observatoire de Paris\newline
77 Av. Denfert Rochereau 75014 PARIS (FRANCE)

\smallskip

\noindent E-mail: \texttt{jacky.cresson@univ-pau.fr}\newline\noindent
\texttt{http://web.univ-pau.fr/\symbol{126}jcresson/}

\bigskip

\noindent Aris DANIILIDIS

\smallskip

\noindent Departament de Matem\`{a}tiques, C1/308\newline Universitat
Aut\`{o}noma de Barcelona\newline E-08193 Bellaterra, Spain.

\smallskip

\noindent E-mail: \texttt{arisd@mat.uab.es} \newline\noindent
\texttt{http://mat.uab.es/\symbol{126}arisd}

\smallskip

\noindent Research supported by the MEC Grant No. MTM2005-08572-C03-03 (Spain).

\bigskip

\noindent Masahiro SHIOTA\smallskip

\noindent Department of Mathematics\newline Nagoya University\thinspace
\ (Furocho, Chikusa) \newline Nagoya 464-8602, Japan.\medskip

\noindent E-mail:\thinspace\texttt{shiota@math.nagoya-u.ac.jp}
\end{document}